\documentclass[12pt, reqno]{amsart}
\usepackage{amsmath, amsthm, amscd, amsfonts, amssymb, graphicx, color}
\usepackage[bookmarksnumbered, colorlinks, plainpages]{hyperref}
\usepackage[T1]{fontenc}
\usepackage[polish, english]{babel}
\usepackage[utf8]{inputenc}

\hypersetup{colorlinks=true,linkcolor=red, anchorcolor=green,
citecolor=cyan, urlcolor=red, filecolor=magenta, pdftoolbar=true}

\usepackage[a4paper,margin=2.5cm]{geometry}

\newtheorem{twr}{Theorem}[section]
\newtheorem{lem}[twr]{Lemma}

\theoremstyle{definition}

\theoremstyle{remark}

\numberwithin{equation}{section}

\DeclareMathOperator{\conv}{conv}
\DeclareMathOperator{\absconv}{absconv}

\DeclareMathOperator{\dist}{dist}

\DeclareMathOperator{\ext}{ext}

\DeclareMathOperator{\re}{Re}

\newcommand\dbm{d_{\mathrm{BM}}}
\usepackage{enumerate}

\begin{document}

\setcounter{page}{1}

\title[A characterization of $\ell_{\infty}^n$ by norm-one projections]{A simple proof of the characterization of $\ell_\infty^n$ by norm-one projections and a stability estimate}
\author[T. Kobos]
{Tomasz Kobos}
\address{Faculty of Mathematics and Computer Science, Jagiellonian University, \L ojasiewicza 6, 30-348 Krak\'ow, Poland}
\email{tomasz.kobos@uj.edu.pl}

\begin{abstract} We give a short and elementary proof of a classical fact that an $n$-dimensional real or complex normed space with the absolute projection constant equal to one is linearly isometric to $\ell_\infty^n$. The proof uses only finite-dimensional convexity and linear algebra, with no infinite-dimensional ingredients. The same approach yields a quantitative stability estimate of order $\varepsilon^{\frac{1}{3}}$ when the absolute projection constant is at most $1+\varepsilon$. This improves Zippin's estimate of order $\varepsilon^{\frac{1}{8}}$ in each fixed dimension as $\varepsilon \to 0$. The real and complex cases are treated by the same arguments. \end{abstract}

\maketitle

\section{Introduction}

All normed spaces considered in this paper are nonzero, and $n$ always denotes a positive integer. For a normed space $X$ over $\mathbb{K}$, where $\mathbb{K}=\mathbb{R}$ or $\mathbb{K}=\mathbb{C}$, and a nonzero linear subspace $Y\subseteq X$, the \emph{relative projection constant} is defined by
$$\lambda(Y,X) = \inf\{\|P\| : \ P:X\to Y \text{ is a bounded linear projection onto } Y\},$$ with the convention that the infimum is $\infty$ if no such projection exists. If $Y$ is finite-dimensional, this infimum is always attained by a certain projection. The \emph{absolute projection constant} $\lambda(X)$ is the supremum of $\lambda(T(X),Z)$ over all normed spaces $Z$ and all linear isometric embeddings $T:X\to Z$. A normed space $X$ is called \emph{$1$-injective} if every bounded linear operator from a subspace of a normed space into $X$ extends to the whole space without increasing its norm. This is known to be equivalent to the requirement that every isometric copy $T(X)\subseteq Z$ of $X$ is the range of a norm-one projection, or, in other words, that $\lambda(T(X),Z)=1$ and the infimum is attained. Consequently, for finite-dimensional spaces $X$, $1$-injectivity of $X$ is equivalent to the equality $\lambda(X)=1$.

The characterization of $1$-injective spaces is a classical result in Banach space theory. By definition, these spaces provide an operator-valued analogue of the Hahn--Banach extension theorem for linear functionals. Whereas norm-preserving extensions need not exist for general operators, they always exist when the range space is $1$-injective. The classical characterization states that a normed space $X$ over $\mathbb{K}$ is $1$-injective if and only if it is linearly isometric to the space of continuous $\mathbb{K}$-valued functions $C(K)$, equipped with the supremum norm, for some extremally disconnected compact Hausdorff space $K$. In the real case, Goodner \cite{goodner} and Nachbin \cite{nachbin} independently proved the characterization under the additional assumption that the unit ball has an extreme point. Kelley \cite{kelley} subsequently removed this assumption, and Hasumi \cite{hasumi} extended Kelley's argument to the complex case.

In the finite-dimensional case we shall consider here, the characterization implies that every $n$-dimensional normed space $X$ satisfying $\lambda(X)=1$ is linearly isometric to $\ell_\infty^n$, that is, $\mathbb{K}^n$ equipped with the maximum norm. The converse follows immediately from the Hahn--Banach theorem. Indeed, suppose that $\ell_\infty^n$ is an isometric subspace of a normed space $Y$. For each $j=1,\ldots,n$, the coordinate functional $f_j: \ell_\infty^n \to \mathbb{K}$ (defined by $f_j(x)=x_j$) has norm one and admits an extension $\widetilde{f}_j:Y\to\mathbb{K}$ of norm one. The required projection $P: Y \to \ell_{\infty}^n$ of norm $1$ is then given by $P(y) = (\widetilde{f}_1(y), \ldots, \widetilde{f}_n(y)) \in \ell_{\infty}^n$.

As we shall see, it suffices to assume that $\lambda(X,Y)=1$ for every $(n+1)$-dimensional normed superspace $Y$ of $X$ to conclude that $X$ is linearly isometric to $\ell_\infty^n$. Thus the characterization can be formulated entirely in finite-dimensional terms, making it natural to seek a direct finite-dimensional proof. The original proofs establish substantially more general results, but their restriction to finite dimensions still involves additional structural machinery. Nachbin \cite{nachbin} first derives the binary intersection property for closed balls and then uses an extreme point of the unit ball to construct a vector-lattice order. The finite-dimensional characterization follows from the classification of finite-dimensional vector lattices. Although this route does not require a representation as a space of continuous functions, it involves several intermediate results; in particular, the passage from the projection property to the binary intersection property uses a minimal-function construction based on Zorn's lemma. Goodner \cite{goodner} also constructs one-dimensional superspaces to obtain geometric restrictions on the unit ball, but subsequently introduces an order structure and applies Clarkson's characterization of real $C(K)$-spaces. Kelley's proof \cite{kelley}, which removes the extreme-point assumption from the general real characterization, uses the weak-star topology on the dual, the Krein--Milman theorem, and embeddings into spaces of continuous functions. These arguments concern real spaces only. Hasumi \cite{hasumi} established the complex analogue by adapting Kelley's method and introducing a continuous-selection theorem.

A different approach is provided by Figiel and Johnson in Section~4 of \cite{figieljohnson}, who obtain a stronger isometric characterization in terms of absolutely summing operators. Their proof uses operator-theoretic machinery such as the Pietsch factorization theorem and trace duality between absolutely summing and integral operator norms. Zippin in Lemma~1 of \cite{zippin1} gives a particularly short convexity proof that the range of a norm-one projection on real $\ell_1^N$ is isometric to $\ell_1^n$, where $n$ is the dimension of the range. By duality, this yields an elementary proof of the characterization for real polyhedral spaces, which admit isometric embeddings into finite-dimensional $\ell_\infty$-spaces. This argument does not directly cover arbitrary finite-dimensional spaces, since such embeddings need not exist.  Zippin also established a substantially stronger quantitative result: a finite-dimensional real space with absolute projection constant sufficiently close to one is close to $\ell_\infty^n$ in the Banach--Mazur distance, with a bound independent of the dimension \cite{zippin1,zippin2} (see also the subsequent erratum \cite{zippin3}). This stability theorem yields the exact characterization by applying it with arbitrarily small positive $\varepsilon$ to a space with projection constant one. Its proof, however, requires a substantial quantitative analysis of projections on finite-dimensional $\ell_1$-spaces, developed across the two papers.

The question of finding a simple direct proof of the finite-dimensional characterization was raised by M. Argerami on MathOverflow in September 2024 \cite{mathoverflow}. In his answer, W.~B.~Johnson observed that there were ``no really simple proofs'' and pointed to several alternative approaches.

Our first goal is to present such a proof. We believe that the argument merits being called simple, in the sense that it is short, direct, and entirely elementary, using only basic finite-dimensional convexity and linear algebra, including the hyperplane separation theorem. Moreover, the same argument covers the real and complex cases simultaneously. As already mentioned, it suffices to consider $(n+1)$-dimensional superspaces of $X$ only. In the real case, this criterion is already implicit in the arguments of Goodner \cite{goodner} and Nachbin \cite{nachbin}. We are not aware of an explicit proof of the corresponding criterion in the complex case.

\begin{twr}
\label{thm:main}
Let $X$ be an $n$-dimensional normed space over $\mathbb{K}$. If the equality $\lambda(X, Y)=1$ holds for every $(n+1)$-dimensional normed superspace $X \subseteq Y$ over $\mathbb{K}$, then $X$ is linearly isometric to $\ell_{\infty}^n$.
\end{twr}

The proof is based on a general construction of isometric embeddings into $(n+1)$-dimensional superspaces, given in Lemma \ref{lem:embed}. It can be viewed as a finite-dimensional adaptation of the canonical embedding into $C(B_{X^*})$: instead of using the whole function space, we retain only the evaluation functions associated with vectors of $X$ and one additional continuous equivariant function. The argument also draws on some other ideas from earlier proofs, including the use of suitably normed one-dimensional extensions and the geometry of the dual unit ball. In the finite-dimensional setting, the relevant description of the dual unit ball is obtained by an elementary convex hull computation. In this sense, it is a synthesis of existing ideas tailored to the finite-dimensional setting. Its main new ingredient is the construction of an embedding in Lemma \ref{lem:embed}, which is used to prescribe scalar values at selected extreme functionals. Extremality then forces a norm-one projection to respect these prescribed values, which are chosen to be incompatible with linearity. The construction thus retains some of the advantages of the function-space embedding without leaving the finite-dimensional setting. For this reason, we believe that the underlying mechanism may have further applications.

The proof of Theorem \ref{thm:main} is presented in Section \ref{section:main}. Our approach can be adapted in a relatively straightforward way to obtain a quantitative stability estimate. In Section \ref{section:stability}, we therefore revisit the stability of $\ell_\infty^n$ as the unique minimizer (up to linear isometry) of the absolute projection constant among $n$-dimensional normed spaces. We recall that the \emph{Banach--Mazur distance} between two $n$-dimensional normed spaces $X$ and $Y$ over $\mathbb{K}$ is defined by
$$\dbm(X, Y) = \inf \|T\|\cdot \|T^{-1}\|,$$
where the infimum is taken over all invertible linear operators $T: X \to Y$. Since we consider only finite-dimensional spaces, this infimum is always attained. Zippin \cite{zippin1,zippin2,zippin3} proved that a finite-dimensional real normed space
whose absolute projection constant is sufficiently close to one is close to
$\ell_\infty^n$, with a bound independent of the dimension. An explicit form of the estimate, attributed to Bernau \cite{bernau}, is given in \cite{gkl} (see the proof of Fact~3.3) and \cite{ntj} (Theorem~33.3): if $X$ is an $n$-dimensional real normed space with $\lambda(X)\le 1+\varepsilon$, where $0\le\varepsilon<17^{-8}$, then
$$\dbm(X,\ell_\infty^n) \leq \frac{1+6\varepsilon^{\frac{1}{8}}} {(1-6\varepsilon^{\frac{1}{8}})(1-17\varepsilon^{\frac{1}{8}})}.$$
In particular, $\dbm(X,\ell_\infty^n) \leq 1+64\varepsilon^{\frac{1}{8}}$ for
$\varepsilon \leq 34^{-8}$. Both the admissible range of $\varepsilon$ and the bound are thus independent of $n$. Our result goes as follows.

\begin{twr}
\label{thm:stability}
Let $X$ be an $n$-dimensional normed space over $\mathbb{K}$ and let $\varepsilon \in \left (0,  (2560n)^{-1} \right ]$. If $\lambda(X, Y)\leq 1 + \varepsilon$ for every $(n+1)$-dimensional normed superspace $X \subseteq Y$ over $\mathbb{K}$, then $$\dbm(X, \ell_{\infty}^n) < 1 + 14(n \varepsilon)^{\frac{1}{3}}.$$
\end{twr}

Thus, in each fixed dimension, our estimate improves the order $\varepsilon^{\frac{1}{8}}$ to $\varepsilon^{\frac{1}{3}}$ as $ \varepsilon \to 0$. However, both the admissible range of $\varepsilon$ and the coefficient in the bound do depend on $n$. Moreover, our hypothesis concerns only $(n+1)$-dimensional superspaces of $X$, and the same argument covers real and complex scalars. The proof again relies on the embedding construction in Lemma \ref{lem:embed}. The key additional ingredient is a quantitative exposed-point lemma (Lemma \ref{lem:functionals}), related to the classical theorem of Straszewicz that exposed points are dense among extreme points of a compact convex set in $\mathbb{K}^n$. Our lemma constructs an exposed functional with a small norming slice that is also almost norming at a prescribed unit vector. To establish the lemma, we choose an auxiliary Euclidean structure using the classical bound $\dbm(X, \ell_2^n) \leq \sqrt{n}$, where $\ell_2^n$ denotes $\mathbb{K}^n$ equipped with the Euclidean norm. Apart from this standard normalization, the proof uses only elementary finite-dimensional arguments.

Throughout the paper, $n$-dimensional normed spaces are always identified with $\mathbb{K}^n$ equipped with a norm. The norm of a space $X$ is denoted by $\|\cdot\|_X$. The usual operator norm is denoted simply by $\|\cdot\|$ when the context is clear. In particular, we write $\|f\|$ for the dual norm of a functional $f\in X^*$. We identify $X^*$ with $\mathbb{K}^n$ using the dual coordinates, i.e., the evaluation is given by $f(x)=\sum_{j=1}^n f_jx_j$ for $x=(x_1, \ldots, x_n) \in X$ and $f=(f_1, \ldots, f_n) \in X^*$. The closed unit ball and the unit sphere of $X$ are denoted by $B_X$ and $S_X$, respectively. For a nonempty set $C\subseteq\mathbb{K}^n$, its \emph{convex hull} and \emph{absolutely convex hull} are defined by
$$\conv C = \left\{ \sum_{j=1}^{N} \alpha_j x_j : N \geq 1, \ x_j \in C, \ \alpha_j \in [0, \infty), \ \sum_{j=1}^{N} \alpha_j = 1 \right\}$$
and
$$\absconv C = \left\{ \sum_{j=1}^{N} \alpha_j x_j : N \geq 1, \ x_j \in C, \ \alpha_j \in \mathbb{K}, \ \sum_{j=1}^{N} |\alpha_j| \leq 1 \right\},$$
respectively. If $C$ is a convex subset of $\mathbb{K}^n$, then a point $x\in C$ is called an \emph{extreme point} of $C$ if the equality $x=ty+(1-t)z$ for $y,z\in C$ and $t\in(0,1)$ implies $x=y=z$. The set of extreme points of $C$ is denoted by $\ext C$. By $\mathbb{T} \subseteq \mathbb{K}$ we denote the set of all scalars of modulus $1$, i.e., the set $\{-1, 1\}$ for $\mathbb{K}=\mathbb{R}$ and the complex unit circle for $\mathbb{K}=\mathbb{C}$.

\section{Proof of Theorem \ref{thm:main}}
\label{section:main}

As mentioned in the introduction, the main idea of the proof of Theorem \ref{thm:main} is to construct suitably normed one-dimensional extensions of $X$. Clearly, for every $x\in X$ we have $\|x\|_X=\max_{g\in B_{X^*}}|g(x)|$. Writing functionals in the standard coordinates of $\mathbb{K}^n$, we add to each $g\in B_{X^*}$ an additional coordinate $\psi(g)$ that depends continuously on $g$. Taking the maximum of the moduli of the resulting functionals $(x,t)\mapsto g(x)+t\psi(g)$ defines a norm on $X\times\mathbb{K}$ under rather mild assumptions on $\psi$. This norm preserves the original norm on $X\times\{0\}$. Crucially, $\psi$ need not be linear, leaving considerable freedom in its choice. In the proof a very specific choice of $\psi$ is used. Geometrically, the construction lifts the dual unit ball into one additional coordinate and takes the convex hull of the resulting graph. Under similar assumptions, the same construction yields $(n+m)$-dimensional superspaces by taking $\psi$ with values in $\mathbb{K}^m$. A map $\psi:B_{X^*}\to\mathbb{K}$ is called \emph{equivariant} if $\psi(\alpha g)=\alpha\psi(g)$ holds for every $\alpha\in\mathbb{T}$ and $g\in B_{X^*}$.

\begin{lem} 
\label{lem:embed} 
Let $X$ be an $n$-dimensional normed space over $\mathbb{K}$, and let $\psi:B_{X^*}\to\mathbb{K}$ be a continuous equivariant map that is not identically zero. Suppose that there exists $r\in(0,1)$ such that $\psi(g)=0$ whenever $\|g\| < r$. Then 
$$\|(x,t)\|_Y = \max_{g\in B_{X^*}}|g(x)+t\psi(g)|$$
defines a norm on $Y=X\times\mathbb{K}$. Moreover, the embedding $X \ni x\mapsto(x,0) \in Y$ is isometric, and, under the identification of $Y^*$ with $X^{*} \times \mathbb{K}$ given by $(g,s)(x,t)=g(x)+st$, we have 
$$B_{Y^*} = \conv \left \{ (g,\psi(g)):g\in B_{X^*} \right \}.$$
\end{lem}

\begin{proof}
We first justify that $\| \cdot \|_Y$ is indeed a norm on $Y$. It is well defined, since the set $B_{X^*}$ is compact and $\psi$ is continuous, so the maximum exists. It is homogeneous and the triangle inequality holds, so it is a seminorm on $Y$. To prove that it is positive for nonzero vectors, suppose first that $x \in X$ is nonzero. We can choose $g \in B_{X^*}$ such that $g(x) \neq 0$ and $\|g\| < r$, so that 
$$\|(x, t)\|_Y \geq |g(x) + t\psi(g) | = |g(x)| > 0.$$
On the other hand, if $x=0$, since $\psi$ is not identically zero, we can choose $g \in B_{X^*}$ such that $\psi(g) \neq 0$. Then $\|(0, t)\|_Y  \geq |t \psi(g)| > 0$ for $t \neq 0$, and hence $\| \cdot \|_Y$ is a norm. Furthermore, it is clear that $\|(x, 0)\|_Y = \|x\|_X$ holds for any $x \in X$, and hence, $X$ embeds isometrically in $Y$.

It remains to prove the representation for $B_{Y^*}$. Let $C:= \{ (g, \psi(g)): \ g \in B_{X^*}\}$. Clearly, $\conv C \subseteq B_{Y^*}$ and it is enough to establish the reverse inclusion. Indeed, suppose, towards a contradiction, that $\varphi \in B_{Y^*} \setminus \conv C$. Since $\psi$ is continuous, $C$ is compact and by the hyperplane separation theorem, there exists $y \in Y$ such that
$$\re(\varphi(y)) > \max_{\gamma \in C}\re(\gamma(y))=\max_{\gamma \in C}|\gamma(y)|= \|y\|_Y \geq |\varphi(y)|,$$
where the first equality follows from the fact that $\alpha C = C$ for any $\alpha \in \mathbb{T}$ by the equivariance of $\psi$. This contradiction shows that $\conv C = B_{Y^*}$.
\end{proof}

\begin{proof}[Proof of Theorem \ref{thm:main}]
Since the dual unit ball $B_{X^*} \subseteq \mathbb{K}^n$ has a nonempty interior and $B_{X^*} = \conv \ext B_{X^*}$ (by the finite-dimensional Krein-Milman theorem), the set $\ext B_{X^*}$ spans $\mathbb{K}^n$. In particular, we can choose $n$ linearly independent functionals $f_1, \ldots, f_n \in \ext B_{X^*}$. If the equality $\absconv \{f_1, \ldots, f_n\} = B_{X^*}$ holds, then
$$\|x\|_X = \max\{|f_j(x)|: \ j=1, \ldots, n\},$$
and consequently, the space $X$ is isometric to $\ell_{\infty}^n$ with the isometry $T: X \to \mathbb{K}^n$ given by $T(x) = (f_1(x), \ldots, f_n(x))$. Otherwise, since $B_{X^*}=\conv\ext B_{X^*}$, there exists $f_{n+1}\in\ext B_{X^*}$ such that $ f_{n+1}\notin\absconv\{f_1,\ldots,f_n\}$. We shall derive a contradiction with this assumption. In particular, for all $j=1, \ldots, n$ we have $\mathbb{T} f_j \cap \mathbb{T} f_{n+1} = \varnothing$. Since the sets $\mathbb{T}f_j$ are compact, we can pick $0 < d < \frac{1}{2}$ such that $\dist(f_j, \mathbb{T} f_{n+1}) > d$ for all $j=1, \ldots, n$ (with the distance being measured using the dual norm of $X$).

For $s \geq 0$ let $u(s) = \max\left \{0, 1 - \displaystyle \frac{s}{d} \right \}$. Thus, $u$ is a continuous function satisfying $u(0)=1$ and $u(s) = 0$ for $s \geq d$. Since $\|f_{n+1}\|=1$, we can choose $x_{n+1} \in S_X$ satisfying $f_{n+1}(x_{n+1})=1$. We define the mapping $\psi: B_{X^*} \to \mathbb{K}$ by
$$\psi(g) = g(x_{n+1})u(\dist(g, \mathbb{T}f_{n+1})).$$
The mapping $\psi$ has the following properties:
\begin{enumerate}
\item $\psi$ is continuous.
\item $\psi$ is equivariant, as $\dist(\alpha g, \mathbb{T}f_{n+1}) = \dist(g, \mathbb{T}f_{n+1})$ holds for every $\alpha \in \mathbb{T}$.
\item $\psi(f_{j})=0$ for every $j=1, \ldots, n$, since $\dist(f_j, \mathbb{T}f_{n+1})>d$. Moreover, $\psi(f_{n+1})=f_{n+1}(x_{n+1})u(0)=1.$
\item $\psi(g) = 0$ for all $g \in B_{X^*}$ with $\|g\| < \frac{1}{2}$. Indeed, for such $g$ we have $\dist(g, \mathbb{T}f_{n+1}) \geq \|f_{n+1}\| - \|g\| > \frac{1}{2} > d$.
\end{enumerate}

Let us now consider the space $Y = X \times \mathbb{K}$ with the norm $\| \cdot \|_Y$ defined in Lemma \ref{lem:embed}. Since $X \times \{0\}$ is a subspace of $Y$ isometric to $X$, by the assumption there exists a linear map $P:Y \to X$ such that $P(x, 0) = x$ for all $x \in X$ and $\|P\|=1$. Fix $j \in \{1, \ldots, n+1\}$ and let $\varphi_j = f_j \circ P \in Y^*.$ Then $\|\varphi_j\| \leq \|f_j\| \cdot \|P\| = 1$, so that $\varphi_j \in B_{Y^*}$. Therefore, by Lemma \ref{lem:embed}, there exist functionals $g_1, \ldots, g_N \in B_{X^*}$ and positive real numbers $\lambda_1, \ldots, \lambda_N$ with sum $1$ such that
$$\varphi_j = \sum_{k=1}^{N} \lambda_k (g_k, \psi(g_k)).$$
Evaluating both sides of this equality on $X \times \{0\}$ yields $f_j = \sum_{k=1}^{N} \lambda_k g_k$. Since $f_j$ is an extreme point of $B_{X^*}$, this gives $g_k = f_j$ for all $k=1, \ldots, N$. Thus, from the fact that $\psi(f_j) = 0$ for $j = 1, \ldots, n$, it now follows that
$$\varphi_j=(f_j, \psi(f_j))=(f_j, 0),$$
for such $j$. Similarly, $\varphi_{n+1}=(f_{n+1}, 1)$. Let us now take $y=P(0,1)\in X$. For $j=1,\ldots,n$, we have 
$$f_j(y)=\varphi_j(0,1)=0.$$
Since $f_1,\ldots,f_n$ were chosen to form a basis of $X^*$, this yields $y=0$. On the other hand, 
$$ f_{n+1}(y)=\varphi_{n+1}(0,1)=1,$$
a contradiction. This completes the proof.

\end{proof}

\section{Stability}
\label{section:stability}

In this section we prove Theorem \ref{thm:stability} concerning the stability of $\ell_\infty^n$. In the proof of Theorem \ref{thm:main}, we used extremality of functionals to force a norm-one projection to respect prescribed scalar values at arbitrarily selected extreme functionals. Turning this argument into a quantitative estimate requires a more precise choice of functionals. Specifically, we seek a sufficiently rich family of functionals admitting small norming slices, from which we can select a well-separated subfamily. We obtain such functionals by a quantitative adaptation of the classical proof of Straszewicz's theorem (see Theorem~18.6 in \cite{rockafellar}), which states that exposed points are dense among the extreme points of a compact convex set in a finite-dimensional space. The main idea is to choose a functional in the dual unit ball whose Euclidean distance from an appropriately chosen point is maximal. This requires an auxiliary Euclidean structure well-comparable to the given norm. For this purpose, we use the classical consequence of the John Ellipsoid Theorem \cite{john}, which gives $\dbm(X, \ell_2^n )\leq \sqrt{n}$ for every $n$-dimensional normed space $X$ over $\mathbb{K}$. Although the John Ellipsoid Theorem and the resulting upper bound on the Banach--Mazur distance are best known in the real setting, they also hold over the complex field (see Theorem~15.3 in \cite{ntj}).

\begin{lem}
\label{lem:functionals}
Let $X$ be an $n$-dimensional normed space over $\mathbb{K}$. Let $x \in X$ be a unit vector and let $L > 0$ be a real number. Then there exist a functional $f \in S_{X^*}$ and a vector $v \in S_X$ such that $f(v)=1$ and the following conditions are satisfied:
\begin{enumerate}[(i)]
\item $f(x) \in \mathbb{R}$ and $f(x) \geq 1 - \frac{n-1}{2L}$.
\item If $g \in B_{X^*}$ and $\re g(v) \geq 1-\delta$ for some $\delta>0$, then $\|f-g\| \leq \sqrt{2\delta(L+n)}$.
\end{enumerate}
\end{lem}

\begin{proof}
We recall that we identify both $X$ and $X^*$ with $\mathbb{K}^n$. For $u\in X$, we denote its norm by $\|u\|_X$, while for $g\in X^*$ we denote its dual norm simply by $\|g\|$. It is a classical consequence of the John Ellipsoid Theorem that $\mathcal{E} \subseteq B_{X^*} \subseteq \sqrt{n}\mathcal{E}$, where $\mathcal{E}$ is the ellipsoid of maximal volume contained in $B_{X^*}$. We choose coordinates on $X^*$, and the corresponding dual coordinates on $X$, such that $\mathcal{E}=\mathbb{B}^n$ is the standard Euclidean unit ball. Then $\|g\|\leq\|g\|_2\leq\sqrt{n}\|g\|$  for $g\in X^*$ and the standard Euclidean norm $\| \cdot \|_2$ on $\mathbb{K}^n$. Let also $[\cdot,\cdot]$ be the standard inner product on $\mathbb{K}^n$, linear in the first variable. Thus  $[g,h]=\sum_{j=1}^n g_j\overline{h_j}$ and $[g,g]=\|g\|_2^2.$ With these coordinate conventions, for every $g\in X^*$ and $u\in X$ we have $g(u)=[g,\overline{u}]$, where $\overline{u}=(\overline{u_1},\ldots,\overline{u_n})$ is the coordinatewise conjugate of $u$.

Let $f \in B_{X^*}$ be a functional with maximal Euclidean distance from the vector $-L\overline{x} \in \mathbb{K}^n$, that is 
$$\|f+L\overline{x}\|_2 = \max_{g \in B_{X^*}} \|g+L\overline{x}\|_2.$$
Thus, for every $g \in B_{X^*}$ we have $\|f+L\overline{x}\|^2_2 \geq \|g+L\overline{x}\|^2_2$. This inequality can be rewritten as
$$\|f\|_2^2 + 2L \re f(x) \geq \|g\|_2^2 + 2L \re g(x).$$
In particular, for $g = \alpha f$, where $\alpha \in \mathbb{T}$ is such that $\alpha f(x)=|f(x)|$ we obtain  $f(x) \in \mathbb{R}$. Moreover, for $g \in S_{X^*}$ satisfying $g(x)=1$ we get
$$n + 2 L f(x) \geq \|f\|_2^2 + 2L f(x) \geq 1 + 2L,$$
which yields (i). 

For (ii), let us note that for any $g \in B_{X^*}$
$$\|f+L\overline{x}\|^2_2 \geq \|g+L\overline{x}\|^2_2 = \|(g-f) + f+L\overline{x}\|^2_2 = \|g-f\|_2^2 + 2\re[g-f, f+L\overline{x}] + \|f + L\overline{x}\|_2^2,$$
so that
$$\re[f-g, f+L\overline{x}] \geq \frac{1}{2} \|f-g\|_2^2.$$
If $w:=\overline{f} + Lx \in \mathbb{K}^n$, then the above inequality becomes
$$\re f(w) \geq \re g(w) + \frac{1}{2} \|f-g\|_2^2,$$
for any $g \in B_{X^*}$. In particular, it is clear that $w \neq 0$, since $w=0$ in the preceding inequality would imply $f=g$ for every $g\in B_{X^*}$, which is impossible. Therefore, for $g_0 \in S_{X^*}$ such that $g_0(w)=\|w\|_X$ we obtain
$$\|w\|_X \geq |f(w)| \geq \re(f(w)) \geq g_0(w) = \|w\|_X,$$
and hence $f(w)=\|w\|_X$. Thus, we have $\|f\|=1$. Let $v: = \frac{w}{\|w\|_X} \in X$. Then $f(v) = 1$ and therefore
$$1 \geq \re g(v) + \frac{1}{2\|w\|_X}\|f-g\|_2^2,$$
for any $g \in B_{X^*}$. In particular, if $\re g(v) \geq 1-\delta,$ then
$$\|f-g\|_2^2 \leq 2\delta\|w\|_X.$$
Since $\|f-g\| \leq \|f-g\|_2$, to conclude the proof it remains to bound the norm $\|w\|_X$ from above. We note that
$$\|w\|_X=f(w) = [f, f + L\overline{x}] = \|f\|^2_2 + L f(x) \leq n \|f\|^2 + L\|f\| = n+L.$$
This completes the proof.

\end{proof}

We are ready to prove Theorem \ref{thm:stability}.

\begin{proof}[Proof of Theorem \ref{thm:stability}] 
We set 
$$\rho=(40n\varepsilon)^{\frac{1}{3}}, \qquad \delta=\frac{\rho^3}{6n}, \qquad L=\frac{n}{2\rho}.$$
From the assumed range of $\varepsilon$ we have 
$$0< \rho \leq \frac{1}{4}, \qquad 0<\delta\leq\frac{1}{384}, \qquad L>\frac n2.$$
For each $x\in S_X$, we choose a functional $f_x\in S_{X^*}$ and a corresponding point $v_x\in S_X$ as in Lemma \ref{lem:functionals}. In particular, $f_x(v_x)=1$ and 
$$f_x(x)=\re f_x(x) \geq 1-\frac{n-1}{2L} \geq 1-\rho.$$
Moreover, 
$$ 2\delta(L+n) = \frac{\rho^2(1+2\rho)}{6} \leq \frac{\rho^2}{4},$$
so the same lemma gives 
\begin{equation}
\label{eq:funct}
\re g(v_x) \geq 1 - \delta \Rightarrow \|f_x-g\| \leq \frac{\rho}{2}, 
\end{equation}
for all $x \in S_X$ and $g \in B_{X^*}$.
Let $F=\{ f_x: \ x\in S_X\} \subseteq S_{X^*}$. Let us choose a subset $\{f_1, \ldots, f_{m} \} \subseteq F$ of the largest possible cardinality that satisfies $\|f_i - \alpha f_j\| > \rho$ for all $1 \leq i < j \leq m$ and $\alpha \in \mathbb{T}$. For each $i$, let $x_i\in S_X$ be such that $f_i=f_{x_i}$, and set $v_i=v_{x_i}$. We observe that for a given $g \in B_{X^*}$, there exists at most one index $i \in \{1, \ldots, m\}$ for which $|g(v_i)| \geq 1 - \delta$. Indeed, let us assume that there are two such indices $i < j$ and choose $\alpha_i, \alpha_j \in \mathbb{T}$ such that $\re(\alpha_ig(v_i))=|g(v_i)|$ and $\re(\alpha_j g(v_j)) = |g(v_j)|$. Then by \eqref{eq:funct} we get $\|f_i - \alpha_i g\| \leq \frac{\rho}{2}$ and $\|f_j - \alpha_j g \| \leq \frac{\rho}{2}$, so that
$$\|f_i - \alpha_i \alpha_j^{-1} f_j \| = \|f_i - \alpha_i g + \alpha_i(g - \alpha_j^{-1} f_j)\| \leq \|f_i - \alpha_i g\| + \|g - \alpha_j^{-1}f_j\| \leq \rho,$$
which contradicts the assumed strict separation by $\rho$. Because $\{f_1, \ldots, f_{m} \}$ is a maximal $\rho$-separated subset, for every $g \in F$ there exists $j \in \{1, \ldots, m\}$ and $\alpha \in \mathbb{T}$ such that $\|g-\alpha f_j\| \leq \rho$. For a given vector $x \in S_X$, applying this for $g:=f_x$ yields
$$|f_j(x)| \geq \re(\alpha f_j(x)) \geq \re(f_x(x)) - \|f_x - \alpha f_j\| \geq 1 -\rho - \rho = 1 - 2 \rho.$$
It follows that
$$(1 - 2\rho) \|x\|_X \leq \max_{j=1, \ldots, m} |f_j(x)| \leq \|x\|_X,$$
for every $x \in X$. Hence, the linear map $T: X \to \ell_{\infty}^m$ given by
$$T(x) = (f_1(x), \ldots, f_m(x)) \in \ell_{\infty}^m$$
is injective, and, consequently, $m \geq n$. If $m=n$, then $T$ is an isomorphism onto $\ell_{\infty}^n$ yielding the estimate
$$\dbm(X, \ell_{\infty}^n) \leq \|T\| \cdot \|T^{-1}\| \leq (1 - 2\rho)^{-1} \leq 1 + 4\rho = 1 + 4(40n \varepsilon)^{\frac{1}{3}} < 1 + 14(n\varepsilon)^{\frac{1}{3}}.$$
Therefore, the result will follow once we establish that $m \leq n$. Suppose, towards a contradiction, that $m\geq n+1$. Since $\dim X^*=n$, the functionals $f_1, \ldots, f_{n+1}$ are linearly dependent. Let $c_1, \ldots, c_{n+1} \in \mathbb{K}$, not all zero, be such that $\sum_{j=1}^{n+1} c_jf_j=0$. For each $j=1, \ldots, n+1$ choose $\gamma_j \in \mathbb{T}$ such that $\gamma_j c_j = |c_j|$ (if $c_j=0$, we choose $\gamma_j$ arbitrarily). We define a function $\sigma: \mathbb{K} \to \mathbb{K}$ by $\sigma(0)=0$ and
$$\sigma(s) = \frac{s}{|s|} \cdot  \min \left  \{ 1, \max\left \{ 0,  \frac{2|s| + 2\delta - 2}{\delta} \right \} \right \} \quad \text{ for } s \neq 0.$$
It is straightforward to verify that $\sigma$ is a continuous function satisfying $\sigma(s)=0$ for $|s| \leq 1- \delta$ and $\sigma(s) = \frac{s}{|s|}$ for $|s| \geq 1 - \frac{\delta}{2}$. Moreover, $|\sigma(s)| \leq 1$ and $\sigma(\alpha s) = \alpha \sigma(s)$ for all $s \in \mathbb{K}$, $\alpha \in \mathbb{T}$. Next, we define a mapping $\psi: B_{X^*} \to \mathbb{K}$ by
$$\psi(g) = \sum_{j=1}^{n+1} \gamma_j\sigma(g(v_j)).$$
The mapping $\psi$ is continuous, equivariant, and satisfies $\psi(g)=0$ for $\|g\| \leq 1 - \delta$ (as then $|g(v_j)| \leq \|g\| \cdot \|v_j\|_X \leq 1 - \delta$ for all $j$). Furthermore, as observed in the previous step, for every $g \in B_{X^*}$, at most one summand in the definition of $\psi(g)$ is nonzero. Since $|\gamma_j|=1$ and $|\sigma(s)|\leq1$, it follows that
\begin{equation}
\label{eq:norma}
|\psi(g)| \leq 1.
\end{equation}
In particular, since $f_j(v_j)=1$ we have $\psi(f_j) = \gamma_j$ for every $j=1, \ldots, n+1$. This shows that $\psi$ is not identically zero and by Lemma \ref{lem:embed} the formula
 $$\|(x,t)\|_Y = \max_{g\in B_{X^*}} |g(x)+t\psi(g)|$$ 
 defines a norm on $Y=X\times\mathbb{K}$. Moreover, the embedding $X \ni x \mapsto (x,0) \in Y$ is isometric, and 
 $$B_{Y^*} = \conv \left \{(g,\psi(g)):g\in B_{X^*} \right \}.$$ 
 By the assumption, there exists a linear map $P:Y\to X$ such that $P(x,0)=x$ for every $x \in X$ and $\|P\|\leq 1+\varepsilon$. Set $y=P(0,1) \in X$. Under the identification $(g,s)(x,t)=g(x)+st$ of $Y^*$ with $X^* \times \mathbb{K}$, we have $f\circ P=(f,f(y))$  for $f\in X^*$. Let us fix $j\in\{1,\ldots,n+1\}$. Since $\|f_j \circ P\|\leq 1 + \varepsilon$, there exist functionals $g_1,\ldots,g_N\in B_{X^*}$ and positive real numbers $\lambda_1,\ldots,\lambda_N$ with sum one such that 
 $$(f_j,f_j(y)) = (1+\varepsilon) \sum_{k=1}^N\lambda_k(g_k,\psi(g_k)).$$ 
 Thus 
 \begin{equation} 
 \label{eq:eval} 
 f_j=(1+\varepsilon)\sum_{k=1}^N\lambda_kg_k \quad \text{ and } \quad  f_j(y)=(1+\varepsilon) \sum_{k=1}^N\lambda_k\psi(g_k). \end{equation} 
 Set $A= \left \{ k: \re g_k(v_j)\geq 1-\frac{\delta}{2} \right \} \subseteq \{1,\ldots,N\}$ and $\beta=\sum_{k\notin A}\lambda_k$. Since $\re(g_k(v_j)) \leq |g_k(v_j)| \leq \|g_k\| \cdot \|v_j\|_X \leq 1$ for all $k=1, \ldots, N$, \eqref{eq:eval} yields
 \begin{align*}
1=\re(f_j(v_j)) &= (1+\varepsilon) \sum_{k=1}^{N} \lambda_k \re(g_k(v_j)) \\
&= (1+\varepsilon)\left (\sum_{k \in A} \lambda_k \re(g_k(v_j)) + \sum_{k \not \in A} \lambda_k \re(g_k(v_j)) \right ) \\
&\leq (1 + \varepsilon)\left( 1-\beta + \left (1 - \frac{\delta}{2} \right )\beta \right )= (1 + \varepsilon) \left ( 1 - \frac{\delta \beta}{2} \right ),
\end{align*}
and hence 
\begin{equation}
\label{eq:beta}
(1+\varepsilon)\beta \leq \frac{2\varepsilon}{\delta}.
\end{equation}
For $k \in A$ we have $|g_k(v_j)| \geq 1 - \frac{\delta}{2} > 1- \delta$. Thus the only nonzero summand in $\psi(g_k)$ is the $j$-th summand, and consequently, $\psi(g_k) = \gamma_j \frac{g_k(v_j)}{|g_k(v_j)|}$. Therefore
$$\re(\overline{\gamma_j} \psi(g_k)) = \frac{\re(g_k(v_j))}{|g_k(v_j)|} \geq 1 - \frac{\delta}{2}.$$
For $k \not \in A$ we can simply estimate $\re(\overline{\gamma_j}\psi(g_k)) \geq -|\psi(g_k)| \geq -1$ by \eqref{eq:norma}. Using \eqref{eq:eval} and \eqref{eq:beta}, we obtain

\begin{align*}
\re(\overline{\gamma_j}f_j(y)) &\geq (1+\varepsilon)\left ( (1 - \beta ) \left (1 - \frac{\delta}{2} \right ) - \beta \right ) = (1 + \varepsilon) \left ( 1 - \frac{\delta}{2} \right ) - (1+\varepsilon)\beta\left( 2 - \frac{\delta}{2} \right ) \\
&\geq (1+\varepsilon) \left ( 1 - \frac{\delta}{2} \right ) - \left( 2 - \frac{\delta}{2} \right ) \frac{2\varepsilon}{\delta} \geq \left ( 1 - \frac{\delta}{2} \right ) - \frac{4\varepsilon}{\delta} \\
&\geq 1 - \frac{\rho^3}{12n} - \frac{24n\varepsilon}{\rho^3} \geq 1 - \frac{1}{768} - \frac{3}{5} > \frac{1}{3}.
\end{align*}
This holds for every $j=1, \ldots, n+1$. Since $c_j = |c_j| \overline{\gamma_j}$, the assumed linear dependence of the functionals $f_j$ yields
$$0 = \re \left ( \sum_{j=1}^{n+1} c_jf_j(y) \right ) = \sum_{j=1}^{n+1} |c_j| \re \left ( \overline{\gamma_j} f_j(y) \right ) > \frac{1}{3} \sum_{j=1}^{n+1} |c_j| > 0,$$
which is the desired contradiction. This completes the proof.
\end{proof}

We do not know the optimal order in this stability problem. Since
$$\frac{\lambda(X)}{\lambda(Y)} \leq \dbm(X,Y)$$
holds for all $n$-dimensional normed spaces $X$ and $Y$ over $\mathbb{K}$ (see Corollary~6 in Section~III.B of \cite{wojtaszczyk}), for every $n$-dimensional $X$ we have $\dbm(X, \ell_{\infty}^n) \geq \lambda(X)$. Thus, if $\lambda(X) = 1 + \varepsilon$, then $\dbm(X, \ell_{\infty}^n) \geq 1 + \varepsilon$, so linear order is the best one could hope for. Moreover, even if a stability estimate of this optimal order could be obtained, it is not clear whether the relevant constants could be chosen independently of $n$, as in Zippin's estimate.

\section*{AI Declaration} 
The author acknowledges substantial assistance from AI tools in developing the arguments and preparing this paper. The author posed the problem of finding a simple proof of the characterization of $\ell_\infty^n$ and a quantitative stability estimate to Claude Opus 5.5 (Anthropic). The mathematical arguments presented in this paper were developed through interaction between the author and Claude Opus 5.5. The author subsequently checked and revised the arguments and rewrote the whole paper in his own words. Microsoft Copilot was used to check and refine the revised arguments, explore further simplifications and assist with language editing. The author independently verified the mathematical arguments and references and takes full responsibility for the content of the publication.

\end{document}